\documentclass[11pt,reqno]{amsart}
\usepackage{latexsym}
\usepackage{multicol}
\usepackage{verbatim}
\usepackage{hyperref}
\usepackage{amsmath, amscd}
\usepackage{color}
\usepackage[latin1]{inputenc}
\usepackage{enumerate}
\usepackage{amsfonts,amsmath,amsthm,amssymb,amscd}
\usepackage{latexsym}
\usepackage{mathrsfs}

\advance\textwidth by 1.2in \advance\oddsidemargin by -.6in \advance\evensidemargin by -.6in
\newtheorem*{cor}{Corollary}
\newtheorem*{lem}{Lemma}

\newtheorem*{prop}{Proposition}

\theoremstyle{definition}

\theoremstyle{definition}
\newtheorem*{thm}{Theorem}

\newtheorem*{rem}{Remark}

\newenvironment{pf}{\proof}{\endproof}
\newcounter{cnt}
\newenvironment{enumerit}{\begin{list}{{\hfill\rm(\roman{cnt})\hfill}}{%
\settowidth{\labelwidth}{{\rm(iv)}}\leftmargin=\labelwidth%
\advance\leftmargin by \labelsep\rightmargin=0pt\usecounter{cnt}}}{\end{list}} \makeatletter
\def\mydggeometry{\makeatletter\dg@YGRID=1\dg@XGRID=20\unitlength=0.003pt\makeatother}
\makeatother \theoremstyle{remark}

\numberwithin{equation}{section}

\let\bwdg\bigwedge
\def\bigwedge{{\textstyle\bwdg}}

\begin{document}

\newcommand{\thmref}[1]{Theorem~\ref{#1}}
\newcommand{\secref}[1]{Section~\ref{#1}}
\newcommand{\lemref}[1]{Lemma~\ref{#1}}
\newcommand{\propref}[1]{Proposition~\ref{#1}}
\newcommand{\corref}[1]{Corollary~\ref{#1}}
\newcommand{\remref}[1]{Remark~\ref{#1}}
\newcommand{\defref}[1]{Definition~\ref{#1}}
\newcommand{\er}[1]{(\ref{#1})}
\newcommand{\id}{\operatorname{id}}
\newcommand{\ord}{\operatorname{\emph{ord}}}
\newcommand{\sgn}{\operatorname{sgn}}
\newcommand{\wt}{\operatorname{wt}}
\newcommand{\tensor}{\otimes}
\newcommand{\from}{\leftarrow}
\newcommand{\nc}{\newcommand}
\newcommand{\rnc}{\renewcommand}
\newcommand{\dist}{\operatorname{dist}}
\newcommand{\qbinom}[2]{\genfrac[]{0pt}0{#1}{#2}}
\nc{\cal}{\mathcal} \nc{\goth}{\mathfrak} \rnc{\bold}{\mathbf}
\renewcommand{\frak}{\mathfrak}
\newcommand{\supp}{\operatorname{supp}}
\newcommand{\Irr}{\operatorname{Irr}}
\renewcommand{\Bbb}{\mathbb}
\nc\bomega{{\mbox{\boldmath $\omega$}}} \nc\bpsi{{\mbox{\boldmath $\Psi$}}}
 \nc\balpha{{\mbox{\boldmath $\alpha$}}}
 \nc\bpi{{\mbox{\boldmath $\pi$}}}
\nc\bsigma{{\mbox{\boldmath $\sigma$}}} \nc\bcN{{\mbox{\boldmath $\cal{N}$}}} \nc\bcm{{\mbox{\boldmath $\cal{M}$}}} \nc\bLambda{{\mbox{\boldmath
$\Lambda$}}}
 \nc\bdd{\rm{ bdd}}

\newcommand{\lie}[1]{\mathfrak{#1}}

\newcommand{\tlie}[1]{\tilde{\mathfrak{#1}}}
\newcommand{\hlie}[1]{\hat{\mathfrak{#1}}}
\newcommand{\tscr}[1]{\tilde{\mathscr{#1}}}
\newcommand{\hscr}[1]{\hat{\mathscr{#1}}}
\newcommand{\hcal}[1]{\hat{\mathcal{#1}}}
\newcommand{\tcal}[1]{\tilde{\mathcal{#1}}}

\makeatletter
\def\section{\def\@secnumfont{\mdseries}\@startsection{section}{1}%
  \z@{.7\linespacing\@plus\linespacing}{.5\linespacing}%
  {\normalfont\scshape\centering}}
\def\subsection{\def\@secnumfont{\bfseries}\@startsection{subsection}{2}%
  {\parindent}{.5\linespacing\@plus.7\linespacing}{-.5em}%
  {\normalfont\bfseries}}
\makeatother
\def\subl#1{\subsection{}\label{#1}}
 \nc{\Hom}{\operatorname{Hom}}
  \nc{\mode}{\operatorname{mod}}
\nc{\End}{\operatorname{End}} \nc{\wh}[1]{\widehat{#1}} \nc{\Ext}{\operatorname{Ext}} \nc{\ch}{\text{ch}} \nc{\ev}{\operatorname{ev}}
\nc{\Ob}{\operatorname{Ob}} \nc{\soc}{\operatorname{soc}} \nc{\rad}{\operatorname{rad}} \nc{\head}{\operatorname{head}}
\def\Im{\operatorname{Im}}
\def\gr{\operatorname{gr}}
\def\mult{\operatorname{mult}}
\def\Max{\operatorname{Max}}
\def\ann{\operatorname{Ann}}
\def\sym{\operatorname{sym}}
\def\Res{\operatorname{\br^\lambda_A}}
\def\und{\underline}
\def\Lietg{$A_k(\lie{g})(\bsigma,r)$}
\def\loc{\operatorname{loc}}

 \nc{\Cal}{\cal} \nc{\Xp}[1]{x(#1)} \nc{\Xm}[1]{y(#1)}
\nc{\on}{\operatorname} \nc{\Z}{{\bold Z}} \nc{\J}{{\cal J}} \nc{\C}{{\bold C}} \nc{\Q}{{\bold Q}}
\renewcommand{\P}{{\cal P}}
\nc{\N}{{\Bbb N}} \nc\boa{\bold a} \nc\bob{\bold b} \nc\boc{\bold c} \nc\bod{\bold d} \nc\boe{\bold e} \nc\bof{\bold f} \nc\bog{\bold g}
\nc\boh{\bold h} \nc\boi{\bold i} \nc\boj{\bold j} \nc\bok{\bold k} \nc\bol{\bold l} \nc\bom{\bold m} \nc\bon{\bold n} \nc\boo{\bold o}
\nc\bop{\bold p} \nc\boq{\bold q} \nc\bor{\bold r} \nc\bos{\bold s} \nc\boT{\bold t} \nc\boF{\bold F} \nc\bou{\bold u} \nc\bov{\bold v}
\nc\bow{\bold w} \nc\boz{\bold z} \nc\boy{\bold y} \nc\ba{\bold A} \nc\bb{\bold B} \nc\bc{\mathbb C} \nc\bd{\bold D} \nc\be{\bold E} \nc\bg{\bold
G} \nc\bh{\bold H} \nc\bi{\bold I} \nc\bj{\bold J} \nc\bk{\bold K} \nc\bl{\bold L} \nc\bm{\bold M} \nc\bn{\mathbb N} \nc\bo{\bold O} \nc\bp{\bold
P} \nc\bq{\bold Q} \nc\br{\bold R} \nc\bs{\bold S} \nc\bt{\bold T} \nc\bu{\bold U} \nc\bv{\bold V} \nc\bw{\bold W} \nc\bz{\mathbb Z} \nc\bx{\bold
x} \nc\KR{\bold{KR}} \nc\rk{\bold{rk}} \nc\het{\text{ht }}

\nc\toa{\tilde a} \nc\tob{\tilde b} \nc\toc{\tilde c} \nc\tod{\tilde d} \nc\toe{\tilde e} \nc\tof{\tilde f} \nc\tog{\tilde g} \nc\toh{\tilde h}
\nc\toi{\tilde i} \nc\toj{\tilde j} \nc\tok{\tilde k} \nc\tol{\tilde l} \nc\tom{\tilde m} \nc\ton{\tilde n} \nc\too{\tilde o} \nc\toq{\tilde q}
\nc\tor{\tilde r} \nc\tos{\tilde s} \nc\toT{\tilde t} \nc\tou{\tilde u} \nc\tov{\tilde v} \nc\tow{\tilde w} \nc\toz{\tilde z}

\def\ch{\operatorname{ch}}
\title{Character formulae and a realization of tilting modules for $\lie{sl}_2[t]$}
\author[ Bennett and Chari  ]{Matthew Bennett and Vyjayanthi Chari}\address{\noindent Department of Mathematics, State University of Campinas, Campinas - SP, 13083-859, Brazil}\address{\noindent Department of Mathematics, University of California, Riverside, CA 92521}
\email{mbenn002@gmail.com}
\thanks{M. B. was partially supported by   FAPESP grant 2012/06923-0 and BEPE Grant 2013/20243-5}
\email{chari@math.ucr.edu}
\thanks{V.C. was partially supported by DMS- 1303052.}
\maketitle
\begin{abstract}
In this paper we study the category of graded modules for the current algebra associated to $\lie{sl}_2$.  The category enjoys many nice properties, including a tilting theory which was established in \cite{BC}.  We  show that  the indecomposable tilting modules for $\lie{sl}_2[t]$  are the  exterior powers of the fundamental global Weyl module and   give the filtration multiplicities in the standard and costandard filtration. An interesting  consequence of our result (which is far from obvious from the abstract definition) is that an indecomposable  tilting module admits a free right action of the ring of symmetric polynomials in finitely many variables. Moreover, if we go modulo the augmentation ideal in this ring, the resulting $\lie{sl}_2[t]$--module is isomorphic to the dual of a  local Weyl module.
\end{abstract}

\section{Introduction}
The current algebra $\lie g[t]$  associated to  a simple Lie algebra $\lie g$ is the algebra of polynomial maps $\bc\to \lie g$. Equivalently,  it is  the complex vector space  $\lie g\otimes\bc[t] $ and the commutator is  the $\bc[t]$--bilinear extension of the Lie bracket of $\lie g$. The Lie  algebra $\lie g[t]$ is graded by the non--negative integers where the $r$--th graded component is $\lie g\otimes t^r$. We are interested in the (non--semisimple)  category $\cal I$  of $\bz$--graded representations of this Lie algebra  with the restrictions  that the graded components are finite--dimensional and,  only finitely many negatively graded components are non--zero.  The  category contains a number of well--known and interesting objects: the $\lie g$--stable Demazure modules occurring in the highest weight integrable representations of the affine Lie algebra associated to $\lie g$ (see \cite{CL}, \cite{FoL}) and also the graded limits of many important families of finite--dimensional representations of quantum affine algebras (see for instance \cite{CM}, \cite{CPweyl},  \cite{Moura}, \cite{Naoi}, \cite{Naoi2}).  Moreover, the category  $\cal I$ is    one of the motivating examples of  affine highest weight categories introduced  recently in \cite{Kh} and \cite{Kl}.

  The isomorphism classes of simple objects in $\cal I$ are indexed by pairs $(\lambda,r)$ where $\lambda$ is a dominant integral weight of $\lie g$ and $r\in\bz$. It is not hard to see that an irreducible module  $V(\lambda, r)$  in the corresponding class is just isomorphic   to the  finite--dimensional irreducible $\lie g$--module associated to $\lambda$. The category is well-behaved in the sense that it has enough projective objects (see \cite{CG})  although these are never finite--dimensional.  The projective cover  $P(\lambda,r)$ of $V(\lambda,r)$  has  two other important quotients: the first one is called the global Weyl module $W(\lambda, r)$ and  is the maximal quotient of $P(\lambda,r)$ whose weights are all bounded above by $\lambda$. The  global Weyl modules are infinite--dimensional  (if $\lambda\ne 0$) and the second quotient which is  of interest to us, is the  unique maximal   quotient of $W(\lambda,r)$ with a one dimensional $\lambda$-weight space.  These are called the local Weyl modules and are denoted as  $W_{\loc}(\lambda,r)$;  they are finite--dimensional,  indecomposable but usually reducible objects of $\cal I$. The definition of the global Weyl module parallels the definition of standard modules which arise in different contexts in the literature  (see for instance \cite{CPS}, \cite{D}).  An important result conjectured in \cite{BCM} and eventually proved in \cite{CI} in complete generality is:  the projective module $P(\lambda,r)$ has a filtration where the successive quotients are global Weyl modules with multiplicities given by the  Jordan--Holder series of the dual  local Weyl modules.   This  implies that the analog of a costandard object in this category should be the dual local Weyl module. (We remark here that the category $\cal I$ is not closed under taking duality,  although  the full subcategory of $\cal I$ consisting of finite--dimensional objects is closed under duality).

With the definitions of standard and costandard objects  in place, it is reasonable to ask if the category contains tilting objects; namely an object which has a filtration where the successive quotients are standard modules and another filtration where the successive quotients are costandard modules.   In \cite{BC}, the authors introduced the notion of an $o$--canonical filtration of an object of $\cal I$ and gave an $\Ext$--vanishing  criterion for the existence of a standard filtration. A similar result for the costandard filtration does not follow since the standard and costandard modules are not dual to each other and in fact such a criterion is not known. This, along with the fact that the index set of simple objects is infinite causes some  difficulty. However,  the main result  of \cite{BC},  shows very abstractly (using the ideas in \cite{Ma} for algebraic groups),  the existence of a unique (up to isomorphism) indecomposable tilting object $T(\lambda, r)$ for $r\in\bz$ and  $\lambda$ a dominant integral weight of $\lie g$.  It was also proved that any other tilting object in $\cal I$ is isomorphic to   the direct sum of  the indecomposable tilting objects.

In this paper we turn to the question of determining the character of the indecomposable  tilting modules in the simplest case of $\lie{sl}_2[t]$. In the case of $\lie{sl}_2$, a dominant integral weight $\lambda$ is just a non--negative integer and we prove that the module $T(\lambda, 0)$  is just a grade shift of   $\wedge^\lambda( W(1,0))$ and we compute its character.  This realization of  the tilting module associated to $(\lambda,0)$  proves that it admits a free  right action  for the ring $\ba_\lambda$  of symmetric polynomials in $\lambda$--indeterminates.  We show  that $$T(\lambda,0)\otimes_{\ba_\lambda} \ba_\lambda/\bi_\lambda\cong_{\lie{sl}_2[t]} W_{\loc}(\lambda, r)^*,$$ for suitable $r\in\bz$, here   $\bi_\lambda$ is  the unique maximal graded ideal of $\ba_\lambda$.
These results should be compared with the  results on global Weyl modules  (see \cite{CPweyl}):  $W(\lambda,0)$  is the  $\lambda$--th  symmetric power of $W(1,0)$ and hence a free  right module for $\ba_\lambda$ with $$W(\lambda,0)\otimes_{\ba_\lambda}\ba_\lambda/\bi_\lambda\cong_{\lie{sl}_2[t]} W_{\loc}(\lambda,0).$$ In the case of higher rank simple Lie  algebras,  suitable analogs of the result for global Weyl modules is known through the work of \cite{CL}, \cite{FoL} and \cite{Naoi3}.

We end the introduction with some comments on the higher rank case.   Based on our result  for $\lie{sl}_2[t]$ it is tempting  to conjecture that  in general,  the tilting module $T(\lambda,0)$ is a free $\ba_\lambda$--module and gives the dual local Weyl module on passing to the quotient by $\bi_\lambda$.
 In the case of $\lie{sl}_{n+1}[t]$ and for a  multiple of the first fundamental weight, it seems likely that our methods will establish this conjecture. In the general case there are several obstructions.  The global Weyl modules $W(\lambda,r)$  have  a  nice presentation which makes the action of the algebra $\ba_\lambda$ on it very natural (see \cite{CFK}) even in the higher rank cases.
However, this is not the case for tilting modules and the abstract construction does not indicate the existence of this action.
   The other difficulty is the following. In the case of $\lie{sl}_2[t]$ we prove in this paper that the tensor product of global Weyl modules has a filtration by global Weyl modules. However, as Mark Shimozono pointed out to us, it is easy to check,  using the known character formulae for the global Weyl module that   the tensor square of the second fundamental representation for $\lie{sl}_{n+1}[t]$ cannot admit a filtration by global Weyl modules.

{\em Acknowledgements.  Part of this work was done when the authors were visiting the Centre de Recherche Mathematique (CRM)  in connection with  the thematic semester \lq\lq New Directions in Lie theory\rq\rq. The authors are grateful to the CRM for the superb working environment.  The authors  are alsovery grateful to the referee for their careful reading of the paper and their many valuable comments and corrections. }

\section{Preliminaries}
\subsection{} Throughout this paper we denote by $\bc$ the field of complex numbers and by $\bz$ (resp. $\bz_+$, $\bn$) the subset of integers (resp. non--negative, positive  integers).  We fix an indeterminate $u$ and  let $\bz[[u,u^{-1}]]$ be the set of power series in $u$ and $u^{-1}$. We shall generally be interested in the case when  the power series is infinite only in one direction.
Given  $n,s\in\bn$, define elements $(1:u)_n\in\bz[[u]]$, and  $[n], \; \qbinom{n}{s}\in\bz[ u]$ by
  \begin{gather*}(1:u)_n=\frac{1}{(1-u)(1-u^2)\cdots (1-u^n)},\\  \\ [n]=\frac{1-u^n}{1-u}, \qquad  \qbinom{n}{s}= \frac{(1-u^n)\cdots (1-u^{n-s+1})}{(1-u)\cdots (1-u^s)}.\end{gather*}
We understand also that $\qbinom{n}{s}=1$ if $s$ is zero and $\qbinom{n}{s}=0$ if $s<0$.

\subsection{} \label{elemgraded}     We say that a complex vector space $V$ is $\bz$--graded  and locally finite--dimensional  if it is  a direct sum of finite--dimensional subspaces $V[r]$, $r\in\bz$.
 Set  $$\mathbb H(V)=\sum_{r\in\bz}\dim V[r]u^r.$$
For $s\in\bz$ we let $\tau^*_s V$ be the grade shift of $V$, i.e. $(\tau^*_s V)[r]= V[r-s]$, it follows that $\mathbb H(\tau_s^*V)= u^{s}\mathbb H(V)$.   The one--dimensional vector space $\bc$ will be regarded  as being $\bz$--graded and located in grade zero.
 We adopt the convention that all unadorned tensor products of complex vector spaces are over $\bc$.   The tensor product of $\bz$--graded vector spaces $V_1$ and $V_2$  is again graded, with grading given by, $$(V_1\otimes V_2)[r]=\bigoplus_{k\in\bz}  V_1[k]\otimes V_2[r-k].$$  If  $\dim V_j[r]=0$ for all $r$ sufficiently large (resp. sufficiently small) for $j=1,2$ then $V_1\otimes V_2$ is locally finite--dimensional  if $V_1$ and $V_2$ are   locally finite--dimensional and in this case we have $$\mathbb H(V_1\otimes V_2)=\mathbb H(V_1)\mathbb H(V_2).$$   A linear map $f: V_1\to V_2$ is graded of degree $p$  if $f(V_1[r])\subset V_2[r+p]$ for some fixed $p\in\bz$. The subspace $\Hom_{\gr}(V_1, V_2) $  of $\Hom_{\mathbb C}(V_1, V_2)$ spanned by  graded   maps  is again a graded vector space. In the special case when $V_1=V$ is locally finite--dimensional, and   $V_2$ is one--dimensional, i.e. $V_2=\tau_s^*\bc$ for some $s\in\bz$, we have $$\mathbb H(\Hom_{\gr}(V,\tau_s^* \bc))=\sum_{r\in\bz}\dim\Hom_{\bc}( V[r],\bc)  u^{s-r}.$$
  {\em In the rest of the paper,  we shall denote by $V^*$ the vector space $\Hom_{\gr}(V,\bc)$ in which  case,  $\tau_s^*(V^*)\cong \Hom_{\gr}(V, \tau_{s}^*\bc)$.}

\subsection{}    Fix an  indeterminate  $t$ and a Lie algebra $\lie a$, and  let $\lie a[t]=\lie a\otimes \bc[t]$ be  the Lie algebra with bracket $$[x\otimes t^r, y\otimes t^s]=[x,y]\otimes t^{r+s},\ \ \ x,y\in\lie a,\ \ r,s\in\bz_+.$$
The Lie algebra $\lie a[t]$ has a natural grading given by the powers of $t$ and the enveloping algebra $\bu(\lie a[t])$ acquires a canonical grading: thus   $\bu(\lie a[t])[k]$ is the subspace of $\bu(\lie a[t])$ spanned by elements of the form $(a_1\otimes t^{r_1})\cdots (a_s\otimes t^{r_s})$ with $s\in\bz_+$ and $r_1+\cdots+r_s=k$. We  identify $\lie a$ with the graded  subalgebra $\lie a \otimes 1 \subset \lie a[t]$.

\subsection{}  By a graded representation of $\lie a[t]$, we mean a $\bz$--graded vector space $V$ such that $$(\lie a\otimes t^r) V[s]\subset V[r+s],\ \ r\in\bz_+,\ \ s\in\bz.$$ If $V_j$, $j=1,2$ are graded representations of $\lie a[t]$, then so are  $V_1\otimes V_2$ and $V_j^*$ with $j=1,2$.  Let   $\ev_0: \lie a[t]\to\lie a$ be   the homomorphism  of Lie algebras  which is the identity on $\lie a$ and zero on $\lie a\otimes t^s$ if $s>0$ .
 The  pull--back of a  representation $V$ of $\lie a$ by  $\ev_0$  gives a graded representation $\ev_0^* V$ of $\lie a[t]$.

\subsection{}\label{gradedprop}  Let $\bc[t_1,\cdots, t_r]$ be the polynomial algebra in $r$--variables and regard it as   being $\bz_+$ graded by requiring the grade of $t_j$ to be one for $1\le j\le r$. The symmetric group  $S_r$ on $r$-letters acts in a natural way on this algebra  by permuting the variables $t_1,\cdots, t_r$.  Let  $\ba_r$  be the subalgebra of symmetric polynomials, and note that $$\mathbb H(\ba_r)=(1:u)_r.$$   Let $\bi_r$ be the maximal ideal of $\ba_r$ defined by taking the sum of all the positively graded components of $\ba_r$. Recall also that $\bc[t_1,\cdots, t_r]$ ia a free $\ba_r$--module of rank $r!$.
{\em In the rest of the paper, we shall use freely the fact  that $\ba_r$ is a polynomial algebra and hence,  by the Quillen--Suslin theorem that a summand of a free module for $\ba_r$ is free.}

\subsection{}  Given a representation $V$ of $\lie a$,  define a representation of $\lie a[t]$ on $V\otimes\bc[t]$, by $$(x\otimes t^r)(v\otimes t^s)= xv\otimes t^{r+s}.$$ It is clear that the obvious grading on  $V\otimes \bc[t]$ makes it  a graded $\lie a[t]$--module. For $r\in\bn$, consider the graded  left $\lie a[t]$--module $(V\otimes\bc[t])^{\otimes r}$.
  Using the identification,  $$(V\otimes\bc[t])^{\otimes r}\cong V^{\otimes r}\otimes \bc[t_1,\cdots, t_r],$$  we see that there exists a  graded   right action of  $\bc[t_1,\cdots, t_r]$  on $(V\otimes\bc[t])^{\otimes r}$ given by right multiplication and hence also a graded  right  action of $\ba_r$. Further $(V\otimes\bc[t])^{\otimes r}$ is a free module for $\ba_r$ under this action.
The action also  commutes  with the left  action of $\lie{a}[t]$, i.e., $(V\otimes\bc[t])^{\otimes r}$ is a graded bimodule for the pair $(\lie a[t], \ba_r)$.

Consider the natural  left action of $S_r$ on  $(V\otimes\bc[t])^{\otimes r}$;  this  commutes with the left action of $\lie a[t]$ and the right action of $\ba_r$.   In particular, every isotypical component of the action of  $S_r$ on $(V\otimes\bc[t])^{\otimes r}$ is a left $\lie a[t]$--module and a right $\ba_r$--module.  Since  $(V\otimes\bc[t])^{\otimes r}$ is a direct sum of the $S_r$--isotypical components it follows that every isotypical component is a free $\ba_r$--module.
We shall be interested in two particular isotypical components, namely the ones  corresponding to the trivial and the sign representation of $S_r$, i.e., the symmetric and wedge powers of $V\otimes \bc[t]$.

\subsection{} In the rest of the paper we take  $\lie a$ to be the Lie algebra $\lie{sl}_2$ of $2\times 2$ complex matrices of trace zero and  $x,y,h$ to be  be the standard basis of $\lie{sl}_2$.  Let $V(1)$ be the natural representation of $\lie{sl}_2$ with basis  $e_j$, $j=1,2$ such that  $$xe_1=ye_2=0,\ \ he_j=(3-2j)e_j, \ \ j=1,2, \ \ ye_1= e_2, \ \ xe_2=e_1.$$ More generally, for  $\lambda\in\bz_+$, let $V(\lambda)$ be the unique (up to isomorphism) irreducible representation of $\lie{sl}_2$ with dimension $\lambda+1$.  It will be convenient to set $V(\lambda)=0$ if $\lambda\notin\bz_+$.

Any finite--dimensional $\lie{sl}_2$--module $V$  is isomorphic to a direct sum of the $V(\lambda)$, $\lambda\in\bz_+$ and $$V\cong\bigoplus_{\mu\in\bz} V_\mu,\ \ V_\mu=\{v\in V: hv =\mu v\}.$$ The character of $V$, denoted as  $\ch V$ is the function $\bz\to\bz_+$ sending $\mu\to\dim V_\mu$. Then,  $$V\cong\oplus_{\lambda\in\bz_+} V(\lambda)^{\oplus m_\lambda}\implies\ \ch V=\sum_{\lambda\in \bz_+}m_\lambda \ch V(\lambda),$$ where, $$ \ch V(\lambda)(\mu) =\begin{cases} 1,\ \ \mu=\lambda-2p,\
 \ 0\le p\le \lambda,\\ 0,\ \ {\rm{otherwise}}.\end{cases}$$

\section{The main results}
\subsection{}  Let $\mathcal I_{\bdd}$ be the category whose objects are  $\bz$--graded representations $V$  of $\lie{sl}_2[t]$  satisfying  $$V[p]=0, \ \ p<<0,\ \ \ \  \dim V[p]<\infty,\ \ p\in\bz,\ \ \ \ V_\mu=0\ \ {\rm{ for\ all\ but\  finitely\ many}}\  \mu\in\bz,$$  and the  morphisms are degree zero maps of $\lie{sl}_2[t]$--modules. The category $\cal I_{\bdd}$ is  abelian and  is closed under tensor products. The dual (see Section \ref{elemgraded}) of a finite--dimensional object of $\cal I_{\bdd}$ is again an object of $\cal I_{\bdd}$.
Given an object $V$ of $\cal I_{\bdd}$, we have $$V\cong_{\lie{sl}_2} \bigoplus_{r\in\bz}  V[r],\qquad \ V_\mu= \bigoplus_{r\in\bz} V[r]_\mu,\ \ \mu\in\bz. $$ The graded character is $$\ \ \ch_{\gr} V=\sum_{r\in\bz} \ch V[r]u^r= \sum_{r\in\bz}\left(\sum_{\lambda\in\bz_+}\dim\Hom_{\lie{sl}_2}(V(\lambda), V[r])\ch V(\lambda)\right)u^r.$$ For brevity, we shall also write $$[V[r]: V(\lambda)]= \dim\Hom_{\lie{sl}_2}(V(\lambda), V[r]).$$
Graded characters are additive on short exact sequences and  multiplicative on tensor products. Moreover, if $V$ is such that  $V^*$ is also an object $\cal I_{\bdd}$, we have
\begin{equation}\label{dualch}\ch_{\gr}V^*= \sum_{r\in\bz}\left(\sum_{\lambda\in\bz_+}\dim\Hom_{\lie{sl}_2}(V(\lambda), V[r])\ch V(\lambda)\right)u^{-r}.\end{equation}  Given an  object of $\cal I_{\bdd}$,  $V$ and $\lambda\in\bz_+$ we let $S^\lambda V$ and $\wedge^\lambda V$ be the $\lambda$--the symmetric and exterior power of $V$ respectively, and note that these are also objects of $\cal I_{\bdd}$.   Recall that the socle (resp. head)  of a module $V$  is the maximal semi--simple submodule, denoted $\soc V$  (resp. quotient $\boh(V$))  of $V$.

\subsection{} We introduce the main objects of interest and summarize the necessary results. We refer the reader to  \cite[Section 2]{CPweyl} (see also \cite{CL}, \cite{FL}) for details. The simple objects of $\cal I_{\bdd}$ are, $$V(\lambda,r)=\tau^*_r\ev_0^* V(\lambda),\ \ r\in\bz,\ \ \lambda\in\bz_+.$$
For  $\lambda\in \bz_+$, the global Weyl module $W(\lambda)$ is defined to be the $\lie{sl}_2[t]$--module generated by an element $w_\lambda $ with relations:\begin{equation}\label{globalweyl}  (x\otimes\bc[t] ) w_\lambda =0,\ \ \  (h\otimes 1)w_\lambda=\lambda w_\lambda,\ \
 \  (y\otimes 1)^{\lambda+1}w_\lambda=0.\end{equation} Since the defining relations of $W(\lambda)$ are graded, it follows that $W(\lambda)$ is a graded ${ \lie{sl}_2[t]}$--module once we define the grade of $w_\lambda$ to be zero.  Moreover, $$W(\lambda)[-r]=0,\ \ r\in\bn, \ \ \dim W(\lambda)[r]<\infty,\ \ W(\lambda)_\mu\ne 0\implies -\lambda\le\mu\le\lambda.$$ It follows that $W(\lambda)\in\cal I_{\bdd}$ and 
 for $r\in\bz$, we set $$W(\lambda,r)= \tau_r^* W(\lambda).$$   The local Weyl module $W_{\loc}(\lambda)$ is the graded quotient of $W(\lambda)$ obtained by imposing the additional relation: $$(h\otimes t^s)w_\lambda=0,\ \ s\in\bn.$$ For $r\in\bz$, we write $W_{\loc}(\lambda,r)=\tau_r^* W_{\loc}(\lambda)$. It is elementary to see that $$\dim W(\lambda, r)[r]_\lambda=\dim W_{\loc}(\lambda,r)_\lambda =1,$$ and we shall denote by $w_\lambda$ any non--zero element of the corresponding space.
 Note that  $\ch_{\gr} W(\lambda,r)=u^r \ch_{\gr} W(\lambda)$ etc.  We use the notation and results of Section \ref{gradedprop} freely in the next proposition.
\begin{prop}\label{sympower}
\begin{enumerit}
\item[(i)] We have an isomorphism of $\lie{sl}_2[t]$---modules: $$W(\lambda)\cong S^\lambda(V(1)\otimes \bc[t]).$$  In particular $W(1)\cong V(1)\otimes\bc[t]$. Moreover,  $W(\lambda)$ is a free right module for $\ba_\lambda$ of rank $2^\lambda$ and $$W_{\loc}(\lambda)\cong W(\lambda)\otimes_{\ba_\lambda} \ba_\lambda/\bi_\lambda.$$
\item[(ii)]  For $0\le p\le \lambda$, we have,  \begin{gather*}\mathbb H(W(\lambda)_{\lambda-2p})=(1:u)_\lambda\mathbb H( W_{\loc}(\lambda)_{\lambda-2p})= (1:u)_\lambda\qbinom{\lambda}{p},\\  \\ \ch_{\gr} W_{\loc}(\lambda)=\sum_{r\in\bz}\left (\qbinom{\lambda}{r}-\qbinom{\lambda}{r-1}\right)\ch V(\lambda-2r),\\ \\
\ch_{\gr} W(\lambda)= (1:u)_\lambda\ch_{\gr} W_{\loc}(\lambda).\end{gather*}
\item[(iii)] The  module $W_{\loc}(\lambda,r)$ has a unique irreducible quotient which is isomorphic to $V(\lambda,r)$.
 Equivalently, the dual module $W_{\loc}(\lambda, -r)^*$ has an irreducible socle which is isomorphic to $V(\lambda, r)$.
 \end{enumerit}\hfill\qedsymbol \end{prop}

\subsection{}  \label{tiltingbasic} We say that an object $M$ of $\cal I_{\bdd}$ has a filtration  by global Weyl modules if there  exists a decreasing  family of submodules $$M= M^0\supset M^1\supset\cdots\supset M^k\supset M^{k+1}= \{0\},$$ such that for all $1\le s\le k$, there exists $\lambda_s\in\bz_+$ and $m_{s,p}\in\bz$, with  $$    M^s/ M^{s+1}\cong \bigoplus_{p\in\bz} W(\lambda_s, p)^{\oplus m_{s,p}}.$$   
It was shown in \cite{BC}, \cite{BCM} that the multiplicity of a global Weyl module   in such a filtration  is independent of the choice of the filtration.
\begin{rem}\label{fininf} It is worth emphasizing that for a fixed $s$,  one could have $m_{s,p}>0$ for infinitely many $p\in\bz_+$. \end{rem}

We shall say that $M$ has a filtration by dual local Weyl modules, if there exists a decreasing   family of submodules $$M=M^0\supset M^1\supset M^2\supset\cdots, \ \ \ \cap_{s\ge 0}M^s= \{0\},$$ such that for all $s\in\bz_+$  there exists  $\lambda_s\in\bz_+$, $p_s\in\bz$. such that $$ \ \ \  M^s/ M^{s+1}\cong  W_{\loc}(\lambda_s, p_s)^*. $$ Again, the multiplicity of the dual local Weyl module  in such a filtration was proved in \cite{BC} to be independent of the filtration. Since the local Weyl modules are finite--dimensional, we see that the filtration is of infinite length unless $\dim M<\infty$.

Say that an object  $T$ of ${\cal I_{\bdd}}$ is tilting if it has a filtration by global Weyl modules and a  filtration by  dual local Weyl modules. Clearly, if $T$ is tilting, then  so is $\tau_s^*T$ for any $s\in\bz$. The dual of the  following was proved in \cite[Theorem 2.7] {BC}.
\begin{thm} \label{tiltexists} Given $\lambda\in\bz_+$ and $r\in\bz$,  there exists a unique (up to isomorphism) indecomposable tilting module $T(\lambda,r)$ such that  \begin{gather} \label{cond1} \wt T(\lambda,r)\subset\{\lambda,\lambda-2,\cdots,-\lambda\},\\  \label{cond2} \dim T(\lambda,r)[r]_\lambda=1, \ \ \ \  \dim T(\lambda,r)[s]_\lambda=0, \ \ s\in\bz, \ \ s<r. \end{gather} In particular,  $T(\lambda,r)=\tau_r^* T(\lambda,0)$. Further,    any tilting module in $\cal I_{\bdd}$ is isomorphic to a direct sum of modules $T(\lambda, r)$,  $\lambda\in \bz_+$, $r\in\bz$.  \hfill\qedsymbol
\end{thm}
\subsection{} The main result of this paper is the following.
\begin{thm} \label{mainthm} Let $\lambda\in\bz_+$ and recall that $W(1)\cong V(1)\otimes\bc[t]$.
We have an isomorphism of objects of $\cal I_{\bdd}$, $$T(\lambda,  a_\lambda)\cong \bigwedge^{\lambda} W(1),\ \ a_\lambda= \binom{\lambda}{2},$$ where we understand that $a_1=0$. In particular,  for all $r\in\bz$, the modules  $T(\lambda,r)$ have  the structure of a free right
 $\ba_\lambda$--module. 
\end{thm}

\subsection{}\label{outline} We outline the main steps of the proof.  Using Theorem \ref{tiltexists} we see that to prove  Theorem \ref{mainthm}, we must show that the $\lie{sl}_2[t]$-module $\wedge^\lambda W(1)$, \\ \\
(1) satisfies  equations \eqref{cond1} and  \eqref{cond2}, with $r$ replaced by $a_\lambda$,\\
(2) is indecomposable,\\
(3) admits a filtration by dual local Weyl modules, and\\
(4) admits a filtration by global Weyl modules.  \\ \\

\subsection{} We  establish  step 1 of the proof. Recall that we have chosen a basis $e_1, e_2$ of $V(1)$.  By Proposition \ref{sympower}(i)  the elements $e_j\otimes t^s$, $j=1,2$, $s\in\bz_+$ are a basis for $W(1)$.  Hence,
 $$W(1)_\mu=0, \ \ \mu\ne \pm 1, \  \ W(1)[s]=0, \ \ s<0,\ \ W(1)[0]_1=\bc.$$
It follows that $$\left( W(1)^{\otimes \lambda}\right) _\mu\ne 0\implies\mu\in\{\lambda,\lambda-2,\cdots,-\lambda\}, $$ which shows that   $\wedge^\lambda W(1)$ satisfies  equation \eqref{cond1}.
Since  $$\left( W(1)^{\otimes \lambda}\right) _{\lambda}= e_1^{\otimes\lambda}\otimes\bc[t_1,\cdots, t_\lambda],$$ it follows that if $f\in\bc[t_1,\cdots, t_\lambda]$, the element $e_1^{\otimes \lambda}\otimes f\in \left(\wedge^\lambda W(1)\right)_\lambda$ iff $f$ is an alternating polynomial in the indeterminates $t_1,\cdots, t_\lambda$. This is equivalent to requiring that $f={\rm{Vand}}(\lambda)g$ for some $g\in\ba_\lambda$, where
if  $\lambda\ge 2$  $, {\rm{Vand}}(\lambda)$  is    the Vandermonde determinant in $t_1,\cdots, t_\lambda$ and ${\rm Vand}(1)=1$. Summarizing, we have shown,  \begin{equation}\label{hilbertwedge}\left(\wedge^\lambda W(1)\right)_\lambda=e_1^{\otimes \lambda}\otimes {\rm Vand}(\lambda)\ba_\lambda, \end{equation}and hence we get $$\left(\wedge^\lambda W(1)\right)[s]_\lambda =\begin{cases} 1,\ \ s=a_\lambda,\\ 0,\ \ s <a_\lambda.\end{cases}$$  This proves that  $\wedge^\lambda W(1)$ satisfies  equations \eqref{cond2} with $r$ replaced by $a_\lambda$.

\subsection{} The proof of steps 2 and 3 can be found in Section \ref{indecomposable}. In the course of establishing these steps we shall also prove the following,
\begin{prop}\label{char} For $\lambda\in P^+$ the module $\wedge^\lambda W(1)$ is a free right $\ba_\lambda$--module, and  we have \begin{gather*} \wedge^\lambda W(1)\otimes_{\ba_\lambda} \ba_\lambda/\bi_\lambda\cong \tau_{a_\lambda}^* W_{\loc}(\lambda)^*,\ \ \ch_{\gr} \wedge^\lambda W(1)= u^{a_\lambda}(1:u)_\lambda \ch_{\gr} W_{\loc}(\lambda)^*.\end{gather*}
\end{prop} 

\subsection{}  The proof of  step 4  is given in Section \ref{global}, where we establish, \begin{prop}\label{motimes1}   If  $M$ has a   filtration by global Weyl modules.  then   $M\otimes W(1)$ also admits such a filtration.  In  particular, $W(1)^{\otimes \lambda}$    has a   filtration by global Weyl modules.
 \end{prop}

\begin{cor} If $M, N\in\Ob  \cal I_{\bdd}$   admit a  filtration by global Weyl modules  then $M\otimes N$ also admits such a filtration.
\end{cor}
\begin{rem} The corollary  is  false for higher rank simple Lie algebras. In fact in the case of $\lie{sl}_4$, it is not hard to check by computing the characters that the tensor square of the global Weyl module associated to the  second fundamental representation  cannot admit a filtration by global Weyl modules.
\end{rem}\noindent  
\subsection{} We isolate the following consequence of  Theorem \ref{mainthm} and Proposition \ref{char}.
\begin{prop} For $\lambda\in \bz_+$ and $r\in\bz$, the tilting module $T(\lambda, r)$ admits a free  right action of $\ba_\lambda$ which commutes with the left action of $\lie{sl}_2[t]$. Moreover, $$\ch_{\gr}T(\lambda, a_\lambda+ r)=u^{a_\lambda}(1:u)_\lambda\ch_{\gr} W_{\loc}(\lambda,-r)^*.$$\hfill\qedsymbol
\end{prop}
\subsection{} Our final result expresses the graded character of the tilting module as a linear combination of the graded character of global Weyl modules.
\begin{lem} \label{final} We have $$\ch_{\gr} T(\lambda, 0
)=\sum_{s=0}^{\lfloor \lambda/2\rfloor} u^{s(s-\lambda)}(1:u)_s\ch_{\gr}W(\lambda-2s).$$
\end{lem} 
The proof of the Lemma is  in Section \ref{globch}.

\section{Indecomposability of $\wedge^\lambda W(1)$ and the dual Weyl filtration}\label{indecomposable}
In this section we prove steps 2 and 3 of Section \ref{outline} and establish Proposition \ref{char}.

\subsection{} The following was proved in \cite{BCGM}    and will play an important role in this section.
\begin{prop} \label{bcgmprop}\hfill
 Suppose that $\mu\in \bz^+$, $r\in\bz$ and that $\psi: \tau_r^* W(\mu)\to W(1)^{\otimes \lambda}$ is any non--zero map. Then $r\ge 0$, $\lambda=\mu$ and $\psi$ is injective.
\hfill\qedsymbol
\end{prop}

\subsection{}  We now prove,
\begin{prop}\label{indec} Let $\lambda\in\bz_+$.
\begin{enumerit}

\item[(i)] If $M$ is a non--zero graded $\lie{sl}_2[t]$--submodule of $\wedge^\lambda W(1)$ then $M_\lambda\ne 0.$
\item[(ii)] The subspace $ (\wedge^\lambda W(1))_\lambda $ is generated as a $(h\otimes \bc[t])$--module by  the element $e_1^{\otimes \lambda}\otimes \rm{Vand}(\lambda)$.
\item[(iii)] The module  $\wedge^\lambda W(1)$ is an indecomposable $\lie{sl}_2[t]$--module.\end{enumerit}\end{prop}
\begin{pf} Suppose that $M$ is a graded $\lie{sl}_2[t]$--submodule of $\wedge^\lambda W(1)$.  Since $\wt M\subset\lambda-\bz_+$ and $M$ is a sum  of finite-dimensional $\lie{sl}_2$-modules,  there exists $\mu\in\bz_+$ and $p\in\bz$, with  $$ M[p]_\mu\ne 0,\ \qquad  (x \otimes  \bc[t])M[p]_\mu=0.$$  It follows from the defininig relations of $W(\lambda)$ (see  equation \eqref{globalweyl}), that  there exists a non--zero map $$\psi: \tau_p^* W(\mu)\longrightarrow M\hookrightarrow\wedge^\lambda W(1)\hookrightarrow W(1)^{\otimes\lambda}.$$ Proposition \ref{bcgmprop} implies that  $\psi$ is injective and $\mu=\lambda$, i.e. $M_\lambda\ne 0$ as required and proves part (i) of the proposition.

Recall that $\ba_\lambda$ is generated as a subalgebra of $\bc[t_1,\cdots, t_\lambda]$  by the  elements $t_1^s+\cdots +t_\lambda^s$, $s\in\bz_+$. Part (ii) follows from \eqref{hilbertwedge} and the fact that $$(h\otimes t^s)\left(e_1^{\otimes \lambda}\otimes{ \rm{Vand}}(\lambda)\right)= e_1^{\otimes \lambda}\otimes{ \rm{Vand}}(\lambda)(t_1^s+\cdots +t_\lambda^s),\ \ \ s\in\bz_+.$$

To prove (iii),  suppose that  $$\wedge^\lambda W(1)= U\oplus V,$$   where $U$ and $V$ are  graded $\lie{sl}_2[t]$-submodules of $\wedge^\lambda W(1)$.  Since  $(\wedge^\lambda W(1)[a_\lambda])_{\lambda}$ is one--dimensional  it must be contained in $U$  or $V$, say $U$. By  part (ii) of the proposition we have  $ (\wedge^\lambda W(1))_\lambda \subset U$ which proves that $V_\lambda=0$.  Part (i) implies that $V=0$ and the proof  is complete.\end{pf}

\subsection{} \label{head}  The local Weyl modules $W_{\loc}(\lambda,r)$ have a natural universal property which is straightforward from the defining relations. Namely, let  $V\in\Ob\cal I_{\bdd}$ be  such that  $$\dim V_\lambda =\dim V[r]_\lambda=p,\ \ \  V_\mu=0,\ \ \mu>\lambda,$$ and assume that  $V$ is generated as an $\lie{sl}_2[t]$--module by $V_\lambda$. Any vector space isomorphism from   $\varphi: W_{\loc}(\lambda, r)_\lambda^{\oplus p}\to V_\lambda$ extends to
 a  surjective morphism of objects of $\cal I_{\bdd}$,  $$\tilde{\varphi}: W_{\loc}(\lambda, r)^{\oplus p}\longrightarrow  V\to 0.\ \ $$ Moreover, the elementary representation theory of $\lie{sl}_2$ implies that we also have  $(\ker\tilde{\varphi})_{-\lambda}=0$.
  We shall frequently need the analog of this discussion for the dual local Weyl module and  we isolate it in the following result.
\begin{lem}\label{dualuniversal} Suppose that $U\in\Ob\cal I_{\bdd}$ is such that $$\dim U_\lambda= \dim U[r]_\lambda =p,\ \  \ \ U_\mu=0,\ \ \mu>\lambda,$$  and assume that every non--zero submodule of $U$ intersects $U_\lambda$ non--trivially. There exists an injective morphism $$0\to U\to \left(W_{\loc}(\lambda, -r)^*\right)^{\oplus p},$$ of objects of $\cal I_{\bdd}$.
\end{lem}
\begin{pf} Let $V$ be the $\lie{sl}_2[t]$-- submodule of $U^*$ generated by $(U^*)_{\pm \lambda}$; note that  the  representation theory of $\lie{sl}_2$  implies that  $(U^*)_{-\lambda}\subset V$
If $V$ is a proper submodule of $U^*$ then the quotient $U^*/V$ is non--zero and $(U^*/V)_\lambda=0$. By taing duals we see that this means that $U$ contains a non--zero submodule which intersects $U_\lambda$ trivially, contradicting our hypothesis. Hence $U^*=V$ and
$\dim(U^*)[-r]_\lambda= \dim( U^*)_\lambda=p$ . We have thus shown that the discussion preceding the Lemma applies to $U^*$ and the lemma follows  by passing to the dual situation.
\end{pf}

\subsection{}\label{socle} We need  additional information about local Weyl modules.
\begin{lem} For $\lambda\in\bz_+$, set $\bar\lambda=0$ if $\lambda$ is even and $\bar\lambda=1$ otherwise. We have $$\soc W_{\loc}(\lambda)\cong V(\bar\lambda,s ),\ \ s=\lfloor\lambda/2\rfloor\lceil\lambda/2\rceil. $$ Equivalently $$\boh(W_{\loc}(\lambda)^*) \cong V(\bar\lambda, -s),\ \  s=\lfloor\lambda/2\rfloor\lceil\lambda/2\rceil.$$
\end{lem}
\begin{pf}  It was proved in  \cite{CL} and  \cite{CPweyl} that the local Weyl module for $\lie{sl}_2[t]$ is isomorphic  to  a Demazure module in a level one representation of the corresponding affine Lie algebra. It is immediate that the socle of the Weyl module is simple and is the $\lie{sl}_2[t]$--submodule generated by the highest weight vector in the level one representation. Since the $\lie{sl}_2$--submodule generated by this highest weight vector is trivial if $\lambda$ is even and two--dimensional if $\lambda$ is odd, it follows that there exists $p\in\bz$ such that  $\soc W_{\loc}(\lambda)\cong V(\bar\lambda, p)$ as $\lie{sl}_2[t]$--modules. To compute the grade in which the socle is concentrated, recall from Proposition \ref{sympower} that $$\sum_{s\in\bz}[W_{\loc}(\lambda)[s]: V(\bar\lambda)]= \qbinom{\lambda}{\lfloor\lambda/2\rfloor}-\qbinom{\lambda}{\lfloor\lambda/2\rfloor-1}.$$ It is straightforward to check that the highest power of $u$ occurring on the right hand side is $\lfloor\lambda/2\rfloor\lceil\lambda/2\rceil$. The sum of all graded components with grade at least  $\lfloor\lambda/2\rfloor\lceil\lambda/2\rceil$ is a submodule of $W_{\loc}(\lambda)$ and hence contains $\soc W_{\loc}(\lambda)$.  It  follows that $\soc W_{\loc}(\lambda)$ is concentrated in grade $\lfloor\lambda/2\rfloor\lceil\lambda/2\rceil$ as required. The second statement of the Lemma is immediate by taking duals.
\end{pf}

\subsection{} We now define the desired filtration.  Set $M=M^0=\wedge^\lambda W(1)$ and  for $s\ge 1$,  let $M^s$ be   the maximal  $\lie{sl}_2[t]$--submodule of $M^{s-1}$ such that $$M^{s}[a_\lambda+s-1]\cap M_\lambda=0.$$  Equivalently, $M^s$, $s\ge 1$  is the maximal submodule of $M_0$ such that $$M^s[a_\lambda+p]\cap M_\lambda=0,\ \  0\le p\le s-1.$$
In particular, we have \begin{equation}\label{wedgefilt} M[a_\lambda+p]_\lambda\subset M^s,\ \ p\ge s,\qquad  \left(M^s/M^{s+1}\right)_\lambda \cong M[a_\lambda+s]_\lambda.\end{equation}
\begin{lem}\label{intzero} We have $$\bigcap_{s \in\bz_+} M^s=\{0\},$$ and hence $\ch_{\gr}M=\sum_{s\in\bz_+}\ch_{\gr}(M^s/M^{s+1})$.
\end{lem} \begin{pf}  Let  $V$ be a submodule of $\cap_s M^s$, in which case, we have by the definition of $M^s$ that   $V[a_\lambda+s]_\lambda = 0$ for all  $s\ge 0$, i.e., $V_\lambda=0$.   Proposition \ref{indec} (i) implies that $V=0$ and hence the desired intersection is zero.  The second statement on characters is standard and is  identical to the proof of \cite[Proposition 2.9(ii)]{BCM}.
 \end{pf}

\subsection{}  We now prove,
\begin{lem}\label{mfilt} \begin{enumerit}
\item[(i)]  For $s\ge 0$, we have an injective morphism $$ 0\to M^s/M^{s+1}\stackrel{\varphi_s}{\longrightarrow}\left(\ \tau^*_{a_\lambda+s}\left (W_{\loc}(\lambda)^*\right)\right)^{\oplus \dim M[a_\lambda+s]_\lambda}.$$
\item[(ii)]  The map $\varphi_0$ gives an isomorphism $$M^0/M^1\cong \tau^*_{a_\lambda}\left(W_{\loc}(\lambda)^*\right).$$
\end{enumerit}
\end{lem}
\begin{pf}  Set $U=M^s/M^{s+1}$,  let $U'$ be a non--zero submodule of $U$ and let $\tilde U'\supsetneq M^{s+1}$  be its preimage in $M^s$.    Equation  \eqref{wedgefilt}  shows that  $\tilde U'\cap M[a_\lambda+s]_\lambda$ has non--zero image in $U$, i.e., $ U'_\lambda\ne 0$. Since $U_\mu=0$ for all $\mu>\lambda$, we   have now  proved that $U$ satisfies the  hypotheses of Lemma \ref{dualuniversal} which establishes the existence of the  morphism $\varphi_s$ with the desired properties.

We now prove part (ii). Using Lemma \ref{socle} we see that  $W_{\loc}(\lambda)^*$  is generated by any  non--zero element of   weight $\bar\lambda$ and  of grade $\left(- \lceil\lambda/2\rceil\lfloor \lambda/2\rfloor\right)$.
Consider the composite map, $$\varphi: M^0\to M^0/M^1\stackrel{\varphi_0}{\longrightarrow } \tau^*_{a_\lambda}\left(W_{\loc}(\lambda)^*\right).$$ By Lemma \ref{socle}  we have  $$\boh\left(\tau^*_{a_\lambda}\left(W_{\loc}(\lambda)^*\right)\right)= V(\bar\lambda, a_\lambda-\lfloor\lambda/2\rfloor\lceil\lambda/2\rceil)= V(\bar\lambda,\lfloor\lambda/2\rfloor\lceil\lambda/2-1\rceil).$$
Hence  to prove that $\varphi_0$ is surjective, it suffices  to prove that $\varphi(\bov)\ne 0$ where
  $$\bov=\sum_{\sigma\in S_\lambda}\sgn(\sigma) \sigma(e_1\otimes e_1t\otimes\cdots\otimes e_1t^{\lceil\lambda/2\rceil-1}\otimes e_2\otimes e_2t\otimes\cdots\otimes e_2 t^{\lfloor \lambda/2\rfloor-1}).  $$ Note that the  weight of $\bov$ is   $\bar\lambda$ and  it is of grade $\left(\lceil\lambda/2-1\rceil\lfloor \lambda/2\rfloor\right)$. A simple calculation shows that  if we set $\bow= (x\otimes t^{\lceil\lambda/2\rceil})^{\lfloor\lambda/2\rfloor}\bov$ then,  $$\bow =\sum_{\sigma\in S_\lambda}\sgn(\sigma)\sigma(e_1\otimes e_1t\otimes\cdots\otimes e_1t^{\lceil\lambda/2\rceil-1}\otimes e_1t^{\lceil\lambda/2\rceil}\otimes e_1t^{\lceil\lambda/2\rceil+1}\otimes\cdots\otimes e_1 t^{\lfloor \lambda/2\rfloor+\lceil\lambda/2\rceil-1}),$$  i.e., $\bow =e_1^{\otimes \lambda}\otimes {\rm Vand}(\lambda)$ .
Since $\bow$ is the unique element (up to scalara) of weight $\lambda$ and grade $a_\lambda$, it follows from \eqref{wedgefilt} that  the image of $\bow$ in $M^0/M^1$ is non--zero.  Since $\varphi_0$ is injective it follows that  $\varphi(\bow)$ is non--zero and hence $\varphi(\bov)\ne 0$ as required.

\end{pf}

\subsection{}  Given two objects $U$ and $V$ of $\cal I_{\bdd}$ we shall say that $\ch_{\gr}U\le\ch_{\gr}V$ if for all $\mu\in\bz_+$ and $r\in\bz$, we have $$[U[r]: V(\mu)]\le[V[r]: V(\mu)].$$
Lemma \ref{mfilt} gives \begin{equation}\label{injective}\ch_{\gr}M^s/M^{s+1}\le u^{a_\lambda+s}\left(\dim  M\left[s+a_\lambda\right]_\lambda\right)\ch_{\gr}W_{\loc}(\lambda)^*.\end{equation} Hence to prove that the filtration $M^s\supset M^{s+1}$, $s\ge 0$ is a filtration by dual local Weyl modules, or equivalently that the  maps $\varphi_s$ are isomorphisms, it suffices to prove that equality holds for sll $s\in\bz_+$ in the previous equation. 
Using \eqref{wedgefilt},  Lemma \ref{intzero},   we get    \begin{gather*}\label{inequality} \ch_{\gr} M=\sum_{s\ge 0}\ch_{\gr}M^s/M^{s+1}\ \le  u^{a_\lambda}\ch_{\gr} W_{\loc}(\lambda)^*\sum_{p\ge 0} \dim  M\left[p+a_\lambda\right]_\lambda u^p.\end{gather*} Since $$ u^{a_\lambda}\sum_{p\ge 0} \dim  M\left[p+a_\lambda\right]_\lambda u^p=\mathbb H(M_\lambda)=u^{a_\lambda}(1:u)_\lambda,$$ where the last equality 
is from \eqref{hilbertwedge}, we get \begin{equation}\label{ineq2}  \ch_{\gr} M\le u^{a_\lambda}(1:u)_\lambda\ch_{\gr} W_{\loc}(\lambda)^*.\end{equation} Clearly, equality holds in \eqref{injective} for all $s\in\bz_+$  iff it holds in \eqref{ineq2} for all $s\in\bz_+$.


\subsection{}\label{surjects}  Recall  that by the discussion in Section \ref{gradedprop} we know that $M=\wedge^\lambda W(1)$ is a free $\ba_\lambda$--module.

 \begin{prop} We have a surjective morphism   of $\lie {sl}_2[t]$--modules, $$\psi: (M\otimes_{\ba_\lambda} \ba_\lambda/\bi_\lambda)\longrightarrow\tau_{a_\lambda}^*(W_{\loc}(\lambda)^*)\to 0,$$ and hence \begin{equation}\label{ge}\ch_{\gr} M\ge u^{a_\lambda} (1:u)_\lambda\ch_{\gr} W_{\loc}(\lambda)^*.\end{equation}
\end{prop}
\begin{pf} Let $\pi: M\to (M\otimes_{\ba_\lambda} \ba_\lambda/\bi_\lambda)\to 0$ be the canonical map of  graded $\lie{sl}_2[t]$--modules.
 Since $M$ is a free $\ba_\lambda$--module we have   that $M_\mu$ is also a free $\ba_\lambda$--module for all $\mu\in \bz$, and $$\rk_{\ba_\lambda} M_\mu=\dim_\bc (M\otimes_{\ba_\lambda} \ba_\lambda/\bi_\lambda)_\mu.$$
 Since $M_\lambda$ is a free $\ba_\lambda$ module of rank one and generated by the element $e_1^{\otimes \lambda}\otimes {\rm{Vand}}(\lambda)$ it follows that $\pi(e_1^{\otimes \lambda}\otimes {\rm{Vand}}(\lambda))\ne 0$ and hence by \eqref{hilbertwedge},  $$\ker\pi\cap M[a_\lambda]_\lambda=0.$$

The definition of $M^1$ implies that $\ker\pi\subset M^1$ and hence we have a non--zero map $M/\ker\pi\to M/ M^1\to 0$. Composing with the map $\varphi_0$ of  Lemma \ref{mfilt}(ii) proves
the existence of the map $\psi$. In particular, it also shows that $$ \ch_{\gr}\left( M\otimes_{\ba_\lambda}\ba_\lambda/\bi_\lambda\right)\ge u^{a_\lambda}\ch_{\gr} W_{\loc}(\lambda)^*.$$ Hence to establish the inequality in \eqref{ge} it suffices to prove that \begin{equation}\label{grcheq}\ch_{\gr} M=(1:u)_\lambda \ch_{\gr}\left( M\otimes_{\ba_\lambda}\ba_\lambda/\bi_\lambda\right).\end{equation}
  The proof of this    is  identical  to the one given in \cite[Sections  3.6-3.7]{BCM} in the context of global Weyl modules, and we recall the argument briefly for the reader's convenience.  For $\mu\in\bz_+$, set  $M_\mu^x=\{v\in M_\mu: xv=0\}.$  Then  by  \cite[Lemma 1.5]{BCM}, we have a decomposition of $M_\mu$ into graded subspaces,  \begin{equation}\label{bcmlemma}\ M_\mu= \left( (y\otimes 1)M\cap M_\mu\right)\oplus M_\mu^x,\ \ {\rm{and\ so}}\ \ \ch_{\gr} M=\sum_{\mu\in\bz_+}\ch V(\mu) \mathbb H(M_\mu^x).\end{equation}
Since the action of $\lie{sl}_2$ commutes with the action of $\ba_\lambda$  it follows that $M_\mu$ and $M_\mu^x$ are graded $\ba_\lambda$-modules. Since  $M$ is a free $\ba_\lambda$--module and hence  $M_\mu$ and $M_\mu^x$ are also free $\ba_\lambda$ modules and so, $$\mathbb H( M_\mu^x)=\mathbb H\left(M^x_\mu\otimes_{\ba_\lambda}\ba_\lambda/\bi_\lambda\right)\mathbb H(\ba_\lambda)=(1:u)_\lambda\mathbb H\left(M^x_\mu\otimes_{\ba_\lambda}\ba_\lambda/\bi_\lambda\right).$$
Substituting in \eqref{bcmlemma}, we get $$\ch_{\gr}M= 
 (1:u)_\lambda \sum_{\mu\in\bz_+}\ch V(\mu)\mathbb H\left(M^x_\mu\otimes_{\ba_\lambda}\ba_\lambda/\bi_\lambda\right)=(1:u)_\lambda\mathbb \ch_{\gr}(M\otimes_{\ba_\lambda}\ba_\lambda/\bi_\lambda),
$$  which establsihes \eqref{grcheq} and completes the proof of the proposition.

\end{pf}

\subsection{} Equation  \eqref{ge}  shows that equality holds in  \eqref{ineq2} and hence also in \eqref{injective}. This completes the proof that $M$ has a filtration by dual local Weyl modules. Using \eqref{grcheq}, we also get $$\ch_{\gr}(M\otimes_{\ba_\lambda}\ba_\lambda/\bi_\lambda)= u^{a_\lambda} \ch_{\gr} W_{\loc}(\lambda)^*.$$ This shows that the map $\psi: (M\otimes_{\ba_\lambda} \ba_\lambda/\bi_\lambda)\longrightarrow\tau_{a_\lambda}^*(W_{\loc}(\lambda)^*)$ is an isomorphism and  we have also proved Proposition \ref{char}.

\section{Proof of Proposition \ref{motimes1} and the global Weyl filtration of $\wedge^\lambda W(1)$.}\label{global}

We begin this section by recalling an important  homological criterion (established in \cite{BC} in the dual case)  for an object of $\cal I_{\bdd}$  to admit a filtration by global Weyl modules. We remark here that  the dual objects lie in the category where the objects have only finitely many positively graded pieces and there is no difficulty in going between these categories using the duality functor.
\subsection{}  Given $M\in\Ob\cal I_{\bdd}$ with $M=\oplus_{\mu\in\bz} M_\mu$, and $\lambda\in\bz_+$,  let $M^\lambda$ be the $\lie {sl}_2[t]$--submodule generated by  the {eigenspaces} $M_{\lambda+\nu}$ with $\nu\in\bz_+$. Then,   $$M^0=M,\ \ \ \ M^\lambda\supset M^{\lambda+1},$$ and $M^\lambda =0$ for $\lambda $ sufficiently large.  This decreasing filtration   is called the $o$--canonical filtration of $M$.
    The dual of the following proposition was established in \cite[Section 3]{BC}.
\begin{prop}\label{oc} Let  $M\in\Ob\cal I_{\bdd}$. For $\lambda\in\bz_+$   with $M^\lambda\ne M^{\lambda+1}$, there exists  a canonical homomorphism of objects of $\cal I_{\bdd}$,$$\bigoplus_{p\in\bz}W(\lambda,p)^{\oplus m_{p,\lambda}}\to M^\lambda/M^{\lambda+1}\to 0, \ \ m_{p,\lambda}\in\bz_+.$$
 Moreover,
 \begin{equation} \label{homeq} m_{p,\lambda}\le \dim\Hom_{\cal I_{\bdd}}(M^\lambda/M^{\lambda+1}, V(\lambda,  p)) = \dim\Hom_{\cal I_{\bdd}}( M, W_{\loc}(\lambda,- p)^*).\end{equation}
and so, we have   \begin{gather}\label{cheq} \ch_{\gr}M=\sum_{\lambda\ge 0}\ch_{\gr}M^{\lambda}/M^{\lambda+1}\le  \sum_{\lambda\ge 0}\sum_{p\in\bz} \dim\Hom_{\cal I_{\bdd}}( M, W_{\loc}(\lambda, -p)^*)\ch_{\gr} W(\lambda, p).\end{gather}\hfill\qedsymbol
\end{prop}

\subsection{} The relation between the $o$--canonical filtration of $M$ and a global Weyl filtration on $M$ was given in  \cite[Propositions  2.6 and 3.11]{BC} and is summarized as follows.
\begin{prop} \label{summand} Let $M\in\Ob{\cal I}_{\bdd}$. The following are equivalent:
\begin{enumerit}
\item[(i)]  $M$ has a filtration by global Weyl modules.
\item[(ii)] The successive quotients in the $o$--canonical filtration are isomorphic to direct sums of global Weyl modules.
and   \begin{gather}\label{cheqo} \ch_{\gr}M=  \sum_{\lambda\ge 0}\sum_{p\in\bz}\dim\Hom_{\cal I_{\bdd}}( M, W_{\loc}(\lambda, -p)^*)\ch_{\gr} W(\lambda, p).\end{gather}
\item[(iii)] For all $(\lambda, p)\in \bz_+\times\bz$, we have $\Ext^1_{\cal I_{\bdd}}(M, W_{\loc}(\lambda, p)^*)=0$.
\end{enumerit}\hfill\qedsymbol
\end{prop}
We note the following corollary and remark that it is non--trivial since the global Weyl  filtrations involved are not finite (see Remark \ref{fininf}) in the sense that $m_{p,\lambda}$ can be non--zero for infinitely many $p$.
\begin{cor}
The  full subcategory consisting  of objects of  $\cal I_{\bdd}$ which admit a filtration by global Weyl modules is closed under
\begin{enumerit}
\item[(i)]   taking direct summands, 
\item[(ii)] extensions and 
\item[(iii)] infinite direct sums.\end{enumerit}
\end{cor}
\begin{pf} All statements follow by applying $\Ext^1_{\cal I_{\bdd}}(-, W_{\loc}(p,r)^*)$ and using Proposition \ref{summand}(iii).
\end{pf}

\subsection{} Assume the following result for the moment.
\begin{prop}\label{wotimesw1} For all $\lambda\in\bz_+$, the module $W(\lambda)\otimes W(1)$ admits a filtration by global Weyl modules.\end{prop}

\subsection{} We complete the proof of  Proposition \ref{motimes1} and its corollary. Let $M$ be an object of $\cal I_{\bdd}$ which admits a filtration by global Weyl modules, say $$M=M^0\supset M^1\supset\cdots\supset
M^k\supset M^{k+1}=\{0\}.$$ It suffices to prove that the successive quotients of  $$M\otimes W(1)\supset M^1\otimes W(1)\supset\cdots\supset M^{k}\otimes W(1)\supset \{0\},$$  are isomorphic to direct sums of global Weyl modules. By Corollary \ref{summand}(ii) it suffices to show that each $M^j\otimes W(1)$ has a filtration by global Weyl modules.
 We prove this by a downward induction on $M^j$.  We fist show that  induction begins at $j=k$ when
$M^k$ is isomorphic to a direct sum of global Weyl modules.  Proposition \ref{wotimesw1} implies   that $M^k\otimes W(1)$ is isomorphic to a direct sum of modules admitting a filtration by global Weyl modules. Using 
 Corollary \ref{summand}(iii) it follows that $M^k\otimes W(1)$
 has a filtration by global Weyl modules.
 Assume the result hold for $M^{j+1}\otimes W(1)$ and
consider the short exact sequence  $$0\to M^{j+1}\otimes W(1) \to M^{j}\otimes W(1) \to \frac{M^j}{M^{j+1}}\otimes W(1)\to 0. $$  Using the induction hypothesis, Corollary \ref{summand}(ii), (iii) and Proposition \ref{wotimesw1} shows that $  M^{j}\otimes W(1)$ has a filtration by global Weyl modules and completes the proof of the inductive step.
 A straightforward induction  on $\lambda$ now  proves that   $W(1)^{\otimes \lambda }$ has a filtration by global Weyl modules.

\medskip

To prove Corollary \ref{motimes1}, an argument identical to the one given above,  shows that it is enough to prove that for all $\lambda,\mu\in\bz_+$ the module $W(\lambda)\otimes W(\mu)$ admits a filtration by global Weyl modules.     Proposition \ref{sympower}(i) shows that $W(\lambda)\otimes W(\mu)$ is a direct summand of $W(1)^{\otimes\lambda+\mu}$ and the result follows by using Corollary \ref{summand}(i).

\subsection{} We complete the proof that $\wedge^\lambda W(1)$ has a filtration by global Weyl modules.  Since $\wedge^\lambda W(1)$ is a  $\lie{sl}_2[t]$--summand of   $W(1)^{\otimes(\lambda}$ it follows again by   Corollary \ref{summand}(i),  that $\wedge^\lambda W(1)$ admits a filtration by global Weyl modules .

The proof  Proposition \ref{wotimesw1} needs  several preliminary results.

\subsection{}  We prove,
\begin{lem}\label{hwlw} Let $\lambda,\mu\in\bz_+$.
\begin{enumerit}\item[(i)]
The subspace $w_\lambda\otimes  W(\mu)_{-\mu}$ generates the $\lie{sl}_2[t]$--module $W(\lambda)\otimes W(\mu)$.
\item[(ii)]  The non-zero quotients  occuring in the $o$--canonical filtration of $W(\lambda)\otimes W(\mu)$ are quotients of direct sum of modules $\tau_s^* W(\lambda+\mu-2p)$ for $s,p \in\bz_+$ and $0\le p\le \min\{\lambda,\mu\}$.\end{enumerit}
\end{lem}
\begin{pf} Since $(x\otimes 1)W(\mu)_\mu=0$,  a  standard $\lie{sl}_2$--argument proves that  $(x\otimes 1)^\mu (y\otimes 1)^\mu w$ is a non--zero multiple of $w$ for all $w\in W(\mu)_\mu$.  Since $\wt W(\mu)\subset\{\mu,\mu-2,\cdots, -\mu\}$ we get  $$W(\mu) = \bu(x\otimes\bc[t])W(\mu)_{-\mu}.$$ Hence $\bu(x\otimes\bc[t])(w_\lambda\otimes W(\mu)_{-\mu}) =w_\lambda \otimes W(\mu)$. Applying elements of $y\otimes\bc[t]$ now proves that  $$W(\lambda)\otimes W(\mu)=\bu(\lie{sl}_2[t])(w_\lambda\otimes W(\mu)_{-\mu}).$$

To prove (ii), set  $M= W(\lambda)\otimes W(\mu)$ and let $M=M^0\supsetneq M^1\supsetneq\cdots M^k\supsetneq M^{k+1}= \{0\}$ be { a re-indexing of} the $o$--canonical filtration so that the successive quotients are non--zero.
 By Proposition \ref{oc},
the modules   $M^p/M^{p+1}$, $0\l p\le k$,  are quotients of the direct sum of modules $W(\lambda_p, r)$, $r\in\bz$, with $\lambda_p\in\lambda+\mu-2\bz_+$  and   $\lambda_p>\lambda_{p-1}$ for all $0\le p\le k$. Assume that $\mu\le \lambda$ without loss of generality.  To complete the proof of part (ii) of the lemma we must show that
$\lambda_0\ge \lambda-\mu$.  By part (i) of the lemma, we see that since $M^0/M^1\ne 0$, the map $M= M^0\to M^0/M^1\to 0$ must be non--zero on the subspace $M_{\lambda-\mu}$ and hence we must have $(M^0/M^1)_{\lambda-\mu}\ne 0$ which proves that $\lambda_0\ge\lambda-\mu$ as required.

\end{pf}
\subsection{}\label{tensorch}   Using Lemma \ref{hwlw} (i)  and the fact that $W(\lambda)\otimes W(1)$ has only non--negatively graded components, it is clear that
\[  \dim \Hom_{\cal I_{\bdd}}(W(\lambda)\otimes W(1), W_{\loc}(\mu,-r)^*)\ne 0\implies\mu=\lambda\pm 1,\; r\in \bz_{\ge 0}.  \]
\begin{lem}  For $\lambda\in\bz_+$, we have
\begin{gather}\label{chseql}  \ch_{\gr} \left(W(\lambda)\otimes W(1)\right)=[\lambda+1]\ch_{\gr}W(\lambda+1)+ (1:u)_1 \ch_{\gr} W(\lambda-1).\end{gather}
\end{lem}
\begin{pf} The proof is a straightforward calculation using the explicit formulae for the graded characters of the global Weyl modules given in Proposition \ref{sympower}(ii) and the fact that the graded character is multiplicative on tensor products.
\end{pf}

\subsection{} \label{gar} Following \cite{Ga}  define elements $P_n\in\bu(h\otimes \bc[t])$, $n\in\bz_+$  recursively,  by: $$P_0= 1,\ \ \ \ P_n=-\frac{1}{n}\sum_{s=0}^{n-1} (h\otimes t^{s+1})P_{n-s-1}.$$  It is easily seen that 
monomials in $P_n$, $n\in\bz_+$ are a basis for $\bu((h\otimes \bc[t])$.
 Recall that the assignment  $\Delta(z)= z\otimes 1+1\otimes z$, $z\in\lie{sl_2}[t]$  defines  the  comultiplication on $\bu(\lie {sl}_2[t])$. A simple induction shows that \begin{equation}\label{comultip} \Delta(P_n)=\sum_{s=0}^n P_s\otimes P_{n-s}.\end{equation}  The following result is a consequence of \cite[Lemma 7.16]{Ga} (see \cite{CPweyl}  for this formulation).
\begin{lem} Suppose that $V$ is a $\lie{sl}_2[t]$--module and let  $v\in V$. Then,
$$(x \otimes\  \bc[t]) v=0\implies (x\otimes t)^s(y\otimes 1)^s v= P_sv,\ \ s\in\bz_+.$$\hfill\qedsymbol
\end{lem}
\subsection{} As a consequence of Lemma \ref{gar}, we see that \begin{equation}\label{fromgar}P_s w_\lambda=0,\ \ s\ge\lambda+1,\end{equation} and hence $W(\lambda)_\lambda$ is spanned by monomials in $\{P_s: 0\le s\le \lambda\}$. Since $P_1=-(h\otimes t)$, we get  $$W(1)[r]= \bc (h\otimes t)^r w_1,
\ \ r\ge 0.$$
\begin{lem}  \label{topgen} For $\lambda\in\bz_+$,  we have \begin{gather*}\left (W(\lambda)\otimes W(1)\right)_{\lambda+1}=\sum_{r\ge 0}\bu(h\otimes\bc[t])(w_\lambda\otimes (h\otimes t^r)w_1).\end{gather*} 

\end{lem}
\begin{pf} It clearly suffices to prove that for all $p\in\bn$, the element $w_\lambda\otimes (h\otimes t)^pw_1$ is in the $(h\otimes\bc[t])$--submodule generated by the elements   $w_\lambda\otimes (h\otimes t^s)w_1$, $0\le s\le \lambda$. Recall that $P_1=-(h\otimes t)$.
 Using \eqref{comultip} and \eqref{fromgar}, we get   \begin{gather*} P_{\lambda+1}(w_\lambda\otimes w_1)= P_\lambda w_\lambda\otimes P_1w_1, \\ P_{\lambda}(w_\lambda\otimes P_1w_1)= P_\lambda w_\lambda\otimes P_1w_1+ P_{\lambda-1}w_\lambda\otimes P_1^2 w_1.\end{gather*}This  proves that $P_{\lambda-1}w_\lambda\otimes P_1^2w_1$  is in the  $\bu((h\otimes \bc[t]))$--module generated by  $w_\lambda\otimes P_1^r w_1$, $r=0,1$. Now using, $$P_{\lambda-1}(w_\lambda\otimes P_1^2w_1)= P_{\lambda-1}w_\lambda\otimes P_1^2 w_1 + P_{\lambda-2}\otimes P_1^3 w_1,$$ we get $P_{\lambda-2}\otimes P_1^3 w_1$ is in the  $\bu((h\otimes \bc[t]))$--module generated by  $w_\lambda\otimes P_1^r w_1$, $r=0,1,2$.  Repeating this argument proves that $w_\lambda\otimes P_1^{\lambda+1}w_1$ is in the  $\bu((h\otimes \bc[t]))$--module generated by  $w_\lambda\otimes P_1^r w_1$, $r=0,1,\cdots,\lambda$. A further repetition of these steps with $w_\lambda\otimes P^sw_\lambda$ for $s\in\bn$  proves the  Lemma. 
\end{pf}
\subsection{}  We now prove,
\begin{lem} \label{dimhom1} For $\lambda\in\bz_+$, we have
\begin{gather*}
 \Hom_{\cal I_{\bdd}}( W(\lambda)\otimes W(1), W_{\loc}(\lambda+1, -r)^*)\ne 0\ \ \implies 0\le r\le \lambda.\end{gather*}
\end{lem}
\begin{pf}
Let  $\varphi: W(\lambda)\otimes W(1)\to W_{\loc}(\lambda+1, -r)^*$ be any non--zero map  of objects of $\cal I_{\bdd}$. Since $$\soc  \left(W_{\loc}(\lambda+1, -r)^*\right)  = V(\lambda+1, r),\ \  \dim ( W_{\loc}(\lambda+1,-r)^*)_{\lambda+1}=1,$$we must have $$\varphi(\left(W(\lambda)_\lambda\otimes W(1)_1\right)[r])\ne 0,\ \qquad  \varphi(\left(W(\lambda)_\lambda\otimes W(1)_1\right)[s])=0,\ \   s\ne r.$$ On the other hand Lemma \ref{topgen} implies that $\varphi(w_\lambda\otimes (h\otimes t)^pw_1)\ne 0$ for some $0\le p\le \lambda$. Hence we must have $p=r$ and the proof is complete.\end{pf}

\subsection{} The next result we need is,
\begin{lem}\label{dimhom2} For $\lambda\in\bz_+$, we have
\begin{gather*} \dim\Hom_{\cal I_{\bdd}}(W(\lambda)\otimes W(1), W_{\loc}(\lambda\pm 1,-r)^*)\le 1,\ \  \ r\ge 0.
\end{gather*}
\end{lem}\begin{pf} To prove the Lemma for $\lambda+1$,  we take $M=W(\lambda)\otimes W(1)$ and note that $$M^{\lambda+1}=\bu(\lie g[t])(W(\lambda)_\lambda\otimes W(1)_1),\qquad \ M^{\lambda+2}=0.$$  Using  equation \eqref{homeq} (with $\lambda$ replaced by $\lambda+1$) we get $$\dim\Hom_{\cal I_{\bdd}}(M, W_{\loc}(\lambda+1,-r)^*)=\dim\Hom_{\cal I_{\bdd}}  (M^{\lambda+1}, V(\lambda+1,r)) ,$$ and hence it is enough to prove that the right hand side is at most one.  If $\varphi: M^{\lambda+1}\to V(\lambda+1, r)$ is any non--zero morphism,  then it is surjective and hence must be non--zero on $W(\lambda)_\lambda\otimes W(1)_1$.
 Lemma \ref{topgen} implies that  the map $\varphi$ is determined by its values on the elements $w_\lambda\otimes (h\otimes t)^pw_1$ for $0\le p\le \lambda$.  If $p\ne r$ the image of $w_\lambda\otimes (h\otimes t)^pw_1$ is zero since $V(\lambda+1,r)$ is concentrated in grade $r$. Hence $\varphi$ is determined by its value on
$w_\lambda\otimes (h\otimes t)^rw_1$ and since  $\dim V(\lambda+1,r)_{\lambda+1}=1 $,  the result follows for $\lambda+1$.

  Let $\varphi: W(\lambda)\otimes W(1)\to W_{\loc} (\lambda- 1, -r)^*$ be any non--zero morphism of $\lie {sl}_2[t]$--modules.  By Lemma \ref{hwlw},  $\varphi$ is determined by its values on $w_\lambda\otimes W(1)_{-1}$. Since $x W(1)_1=0$, and $W(1)$ is isomorphic to the direct sum of finite--dimensional $\lie{sl}_2$--modules, it follows that $y: W(1)_1\to W(1)_{-1}$ is an isomorphism of graded spaces and, hence $$\dim W(1)[s]_{-1}= \dim \bc (y(h\otimes t)^s) w_\lambda\le 1.$$
Since $\dim\left( W_{\loc}(\lambda-1, r)\right)^*[s]_{\lambda-1}\le 1$ with equality holding iff  $s=r$,
we see  that
$$\varphi(w_\lambda\otimes W(1)_{-1})=\varphi(w_\lambda\otimes y(h\otimes t^r)w_1) \subset \left(W_{\loc}(\lambda-1,-r)^*\right)[r]_{\lambda-1},$$ which completes the proof of the Lemma. \end{pf}
\subsection{}{\em Proof  of Proposition \ref{wotimesw1}}. \ \ 
 To prove that $W(\lambda)\otimes W(1)$ admits a filtration by global Weyl modules, we  use Lemma  \eqref{dimhom1} and  Lemma \eqref{dimhom2} and  Proposition \ref{oc}, to get
\begin{gather*}\ch_{\gr}\left(W(\lambda)\otimes W(1)\right)\le
 \sum_{r=0}^\lambda\ch_{\gr} W(\lambda+1, r) +\sum_{r\ge 0}\ch_{\gr} W(\lambda-1, r)\end{gather*}i.e., $$ \ch_{\gr}\left(W(\lambda)\otimes W(1)\right)\le\left(\sum_{r=0}^\lambda u^r\right)\ch_{\gr} W(\lambda+1)+\left(\sum_{r\ge 0}u^r\right)\ch_{\gr}W(\lambda-1).$$ Equation  \ref{chseql} implies that we must have  an equality. Proposition \ref{summand} now gives that   that the $o$--canonical filtration is a  filtration by global Weyl modules,  more precisely, we have the following short exact sequence of objects of $\cal I_{\bdd}$:
$$0\to\bigoplus_{r=0}^\lambda  W(\lambda+1,r)\to W(\lambda)\otimes W(1) \to \bigoplus_{r\ge 0} W(\lambda-1, r)\to 0, $$ proving that  $W(\lambda)\otimes W(1)$ admits a filtration by global Weyl modules.

\subsection{} We determine explicitly the filtration multiplicities in the tensor product of global Weyl modules.

\begin{prop}  For $\lambda,\mu,\nu\in\bz_+$, the multiplicity of $W(\lambda+\mu-2\nu)$ in a global Weyl filtration of $W(\lambda)\otimes W(\mu)$ is $\qbinom{\lambda+\mu-2\nu}{\mu-\nu}(1:u)_\nu,$ or equivalently,
\begin{gather} \ch_{\gr}\left( W(\lambda)\otimes W(\mu) \right)= \sum_{\nu=0}^{\min\{\lambda,\mu\}} \qbinom{\lambda+\mu-2\nu}{\mu-\nu}(1:u)_\nu \ch_{\gr} W(\lambda+\mu-2\nu). \end{gather}
\end{prop}
\begin{pf} Notice first that since we know that $W(\lambda)\otimes W(\mu)$ has a filtration by global Weyl modules, the two statements in the proposition are indeed equivalent. Assume without loss of generality that $\lambda\ge \mu$ and note that the case when $\mu=1$ was  proved in Lemma \ref{dimhom2}.  We complete the proof by induction on $\mu$ and have only to establish the inductive step.  Assuming the result for $\mu$ we prove the result for $\mu+1$ by using  the induction hypothesis to compute the graded character of  $W(\lambda)\otimes W(\mu)\otimes W(1)$ in two ways.
Thus, we have
\begin{gather*}
\ch_{\gr}\left((W(\lambda)\otimes W(\mu))\otimes W(1)\right)=\left(\sum_{\nu=0}^{\mu}\qbinom{\lambda+\mu-2\nu}{\mu-\nu}(1:u)_\nu\ch_{\gr} W(\lambda+\mu-2\nu)\right)\ch_{\gr}W(1)\\
=\left(\sum_{\nu=0}^{\mu}\qbinom{\lambda+\mu-2\nu}{\mu-\nu}(1:u)_\nu[\lambda+\mu+1-2\nu]\ch_{\gr} W(\lambda+\mu+1-2\nu)\right)\\
 +\left(\sum_{\nu=0}^{\mu}\qbinom{\lambda+\mu-2\nu}{\mu-\nu}(1:u)_\nu(1:u)_1\ch_{\gr} W(\lambda+\mu-1-2\nu)\right),
\end{gather*}
and also
\begin{gather*} \ch_{\gr}(W(\lambda)\otimes(W(\mu)\otimes W(1)))= \ch_{\gr} W(\lambda)([\mu+1]\ch_{\gr}W(\mu+1)+ (1:u)_1\ch_{\gr} W(\mu-1)),\\
=[\mu+1]\ch_{\gr}(W(\lambda)\otimes W(\mu+1)) + (1:u)_1\sum_{\nu=0}^{\mu-1}\qbinom{\lambda+\mu-1-2\mu}{\mu-1-\nu}(1:u)_\nu \ch_{\gr} W(\lambda+\mu-1-2\nu).
\end{gather*} Equating the two solutions gives the result.
{  This is easily seen by using the following identity, which is straightforward to establish
\[ (1:u)_\nu(1:u)_1 \left( \qbinom{\lambda+\mu -2\nu}{\mu-\nu}- \qbinom{\lambda+\mu -1-2\nu}{\mu-1-\nu} \right)\]
\[ = (1:u)_{\nu+1} \left( \qbinom{\lambda + \mu -1 -2\nu}{\mu-\nu} [\mu+1] - \qbinom{\lambda + \mu -2-2\nu}{\mu -1-\nu} [\lambda + \mu -1 -2\nu] \right) \].
}
\end{pf}

\subsection{}\label{globch} We now prove Lemma \ref{final}. Using the relation between the graded character of the global and local Weyl modules given in Proposition \ref{sympower}(ii) and Proposition \ref{char}  we see that proving Lemma \ref{final} is  equivalent to proving that $$\ch_{\gr} W_{\loc}(\lambda)^*=\sum_{s=0}^{\lfloor\lambda/2\rfloor}u^{s(s-\lambda)}\frac{(1:u)_s(1_u)_{\lambda-2s}}{(1:u)_\lambda}\ch_{\gr} W_{\loc}(\lambda-2s).$$
For,  $\lambda\in\bz_+$ we see that  Proposition \ref{sympower}(ii) and Lemma \ref{tensorch} imply, \begin{equation}\label{tpweyl}\ch_{\gr} W_{\loc}(\lambda)\ch_{\gr} W_{\loc}(1)=\ch_{\gr} W_{\loc}(\lambda+1)+(1-u^\lambda)\ch_{\gr} W_{\loc}(\lambda-1). \end{equation} Noting that $W_{\loc}(1)$ is self--dual, we get 
 \begin{equation}\label{tpdual}\ch_{\gr} W_{\loc}(\lambda)^*\ch_{\gr} W_{\loc}(1)^*=\ch_{\gr} W_{\loc}(\lambda+1)^*+(1-u^{-\lambda})\ch_{\gr} W_{\loc}(\lambda-1)^*.\end{equation}
Since $$\ch_{\gr} W_{\loc}(\lambda)= \ch V(\lambda)\oplus_{0\le \mu<\lambda} f_\mu\ch V(\mu),$$  for a unique choice of  $f_\mu\in\bz[u]$, it follows that we can write $\ch V(\lambda)-\ch_{\gr} W_{\loc}(\lambda)$ uniquely  as a $\mathbb Z[u]$--linear combination of $\ch_{\gr}W_{\loc}(\mu)$ with $0\le \mu< \lambda$. It follows therefore  that for all $\lambda\in\bz_+$  we can write   \begin{equation}\label{dualloc}\ch_{\gr} W_{\loc}(\lambda)^*=\sum_{s\ge 0} b_\lambda(s) \ch_{\gr} W_{\loc}(\lambda-2s),\end{equation}for a unique choice of $b_\lambda(s)\in\bz[q,q^{-1}]$,  with   $ b_\lambda(0)=1$ and $b_\lambda(s)=0$ if $2s>\lambda$.  Substituting  in equation \eqref{tpdual},  we get \begin{eqnarray*}\sum_{s\ge 0} b_\lambda(s) \ch_{\gr} W_{\loc}(\lambda-2s)\ch_{\gr} W_{\loc}(1)&=&\sum_{s\ge 0}b_{\lambda+1}(s)\ch_{\gr}W_{\loc}(\lambda+1-2s)\\ &+& (1-u^{-\lambda})\sum_{s\ge 0}b_{\lambda-1}(s)\ch_{\gr}W_{\loc}(\lambda-1-2s).\end{eqnarray*} Using equation \ref{tpweyl} on the left hand side of the preceding equation, gives\begin{gather*}\sum_{s\ge 0} b_\lambda(s)\ch_{\gr} W_{\loc}(\lambda+1-2s)+ \sum_{s\ge 0} (1-u^{\lambda-2s})b_\lambda(s)\ch_{\gr}W_{\loc}(\lambda-2s-1)=\\\sum_{s\ge 0}b_{\lambda+1}(s)\ch_{\gr}W_{\loc}(\lambda+1-2s)+ (1-u^{-\lambda})\sum_{s\ge 0}b_{\lambda-1}(s)\ch_{\gr}W_{\loc}(\lambda-1-2s)\end{gather*}
 Equating the  coefficients of $\ch_{\gr} W_{\loc}(\lambda+1-2s)$ on both sides, gives the following recursion: $$b_{\lambda+1}(s)+(1-u^{-\lambda})b_{\lambda-1}(s-1)= b_\lambda(s)+(1-u^{\lambda-2s+2})b_{\lambda}(s-1),$$ with initial conditions $b_0(0)=b_1(0)=1$ and $b_\lambda(s)=0$, $2s>\lambda$.  It is straightforward to check that taking $$b_\lambda(s)= u^{s(s-\lambda)}\frac{(1:u)_s(1:u)_{\lambda-2s}}{(1:u)_\lambda}$$ satisfies the recursion and has the same initial conditions. This proves the Lemma.

\end{document}